\documentclass[aap,seceqn,dvips]{arximspdf}
\usepackage{stmaryrd}
\usepackage{mathbh}

\doi{10.1214/105051606000000222}
\volume{16}
\issue{3}
\pubyear{2006}
\firstpage{1633}
\lastpage{1652}

\makeatletter
\newcommand{\eqref}[1]{\textup{(\ref{#1})}}

\newtheorem{lemma}{Lemma}[section]

\newtheorem{corollary}{Corollary}
\newproclaim{remark}{Remark}

\newtheorem{theorem}{Theorem}[section]

\newproclaim{example}{Example}[section]
\newproclaim{ass}{Assumption}
\makeatother

\begin{document}
\begin{frontmatter}

\title{A filtering approach to tracking volatility from prices observed at random times}
\runtitle{Filtering random-times observations}

\begin{aug}
\author[A]{\fnms{Jak\v sa} \snm{Cvitani\'c},\corref{}\ead[label=e1]{cvitanic@hss.caltech.edu}\thanksref{t1}}
\author[B]{\fnms{Robert} \snm{Liptser}\ead[label=e2]{liptser@eng.tau.ac.il}} and
\author[C]{\fnms{Boris} \snm{Rozovskii}\ead[label=e3]{rozovski@math.usc.edu}\thanksref{t2}}
\pdfauthor{J. Cvitanic, R. Liptser, B. Rozovskii}
\runauthor{J. Cvitani\'c, R. Liptser and B. Rozovskii}
\thankstext{t1}{Supported in part by the NSF Grants DMS-00-99549 and DMS-04-03575.}
\thankstext{t2}{Supported in part by the Army Research
Office and the Office of Naval Research under Grants DAAD19-02-1-0374 and N0014-03-0027.}
\affiliation{Caltech, Tel Aviv University and University of Southern California}
\address[A]{J. Cvitani\'c\\
Caltech, M/C 228-77\\
1200 E. California Blvd.\\
Pasadena, California 91125\\
USA\\
\printead{e1}}
\address[B]{R. Liptser\\
Department of Electrical Engineering-Systems\\
Tel Aviv University\\
69978 Tel Aviv\\
Israel\\
\printead{e2}}
\address[C]{B. Rozovskii\\
Department of Mathematics\\
University of Southern California\\
Los Angeles, California 90089-1113\\
USA\\
\printead{e3}}
\end{aug}

\received{\smonth{12} \syear{2003}}
\revised{\smonth{2} \syear{2006}}

\begin{abstract}
This paper is concerned with nonlinear filtering of the coefficients in asset
price models with stochastic volatility. More specifically, we assume that the
asset price process $S=(S_{t})_{t\geq0}$ is given by
\[
dS_{t}=m(\theta_{t})S_{t}\,dt+v(\theta_{t})S_{t}\,dB_{t},
\]
where $B=(B_{t})_{t\geq0}$ is a Brownian motion, $v$ is a positive function
and $\theta=(\theta_{t})_{t\geq0}$ is a c\'{a}dl\'{a}g strong Markov process.
The random process $\theta$ is unobservable. We assume also that the asset
price $S_{t}$ is observed only at random times $0<\tau_{1}<\tau_{2}<\cdots.$
This is an appropriate assumption when modeling high frequency financial data
(e.g., tick-by-tick stock prices).

In the above setting the problem of estimation of $\theta$ can be
approached as a special nonlinear filtering problem with measurements
generated by a multivariate point process $(\tau_{k},\log
S_{\tau_{k}})$. While quite natural, this problem does not fit into the
``standard'' diffusion or simple point process filtering frameworks and
requires more technical tools. We derive a closed form optimal
recursive Bayesian filter for $\theta_{t}$, based on the observations
of $(\tau_{k},\log S_{\tau_{k}})_{k\geq1}$. It turns out that the
filter is given by a recursive system that involves only deterministic
Kolmogorov-type equations, which should make the numerical
implementation relatively easy.
\end{abstract}


\end{frontmatter}

\section{Introduction}\label{Sec-1}

In the classical Black--Scholes model for financial markets, the stock
price $S_{t}$ is modeled as a geometric Brownian motion, that is, with
diffusion coefficient equal to $\sigma S_{t}$, where ``volatility''
$\sigma$ is assumed to be constant. The volatility parameter is of
great importance in applications of the model, for example, for option
pricing. Consequently, many researchers have generalized the constant
volatility model to so-called stochastic volatility models, where
$\sigma_{t}$ is itself random and time dependent. There are two basic
classes of models: complete and incomplete. In complete models, the
volatility is assumed to be a functional of the stock price; in
incomplete models, it is driven by some other source of noise that is
possibly correlated with the original Brownian motion. In this paper we
study a particular incomplete model in which the volatility process is
independent of the driving Brownian motion process. This has the
economic interpretation of the volatility being influenced by market,
political, financial and other factors that are independent of the
``systematic risk'' (the Brownian motion process) associated with the
particular stock price under study.
Option traders, investment banks, economic analysts and others depend
on modeling future volatility for their trading, economic forecasts,
risk management and so on.

Estimating volatility from observed stock prices is not a trivial task
in either complete or incomplete models, in part because the prices are
observed at discrete, possibly random time points. Since volatility
itself is not observed, it is natural to apply filtering methods to
estimate the volatility process from historical stock price
observations. Nevertheless, this has only recently been investigated in
continuous-time models, in particular, by Frey and Runggaldier
\cite{FR}. See~\cite{R} for an up-to-date survey. See also \cite{EHJ}
for a discrete-time approach with equally spaced observations,
\cite{GT} for an approximating algorithm in continuous time, \cite{MM}
for a nonparametric approach, as well as \cite{FPS,KX,RZ} for still
other approaches. There is also a rich econometrics, time-series
literature on ARCH--GARCH models of stochastic volatility, that
presents an alternative way to model and estimate volatility;
see~\cite{Gou} for a survey.

Our paper was prompted by Frey and Runggaldier \cite{FR}. Like that paper, we
assume that the asset price process $S=(S_{t})_{t\geq0}$ is given by
\[
dS_{t}=m(\theta_{t})S_{t} \,dt+v(\theta_{t})S_{t}\,dB_{t},
\]
where $B=(B_{t})_{t\geq0}$ is a Brownian motion, $v$ is a positive
function, and $\theta=(\theta_{t})_{t\geq0}$ is a c\'{a}dl\'{a}g strong
Markov process. The ``volatility'' process $\theta$ is unobservable,
while the asset price $S_{t}$ is observed only at random times
$0<\tau_{1}<\tau_{2}<\cdots.$ This assumption is designed to reflect the
discrete nature of high frequency financial data such as tick-by-tick
stock prices. The random time moments $\tau_{k}$ can be interpreted as
``instances at which a large trade occurs or at which a market maker
updates his quotes in reaction to new information'' (see \cite{F}).
Hence, it is natural to assume that $\{ \tau_{k}\} _{k\geq1}$ might
also be correlated with $\theta.$

In the above setting the problem of volatility estimation can be regarded as a
special nonlinear filtering problem.

Frey and Runggaldier \cite{FR} derive a Kallianpur--Striebel type
formula (see, e.g.,~\cite{KalStr}) for the optimal mean-square filter
for $\theta_{t}$ based on the observations of
$S_{\tau_{1}},S_{\tau_{2}},\ldots$ for all $\tau_{k}\leq t$ and
investigate Markov chain approximations for this formula. We extend
this result in that we derive the exact filtering equations for
$\theta_{t}$ that allow us to compute the conditional distribution of
$\theta_{t}$ given $S_{\tau_{1}\wedge t}$, $S_{\tau_{2}\wedge
t},\ldots.$ Moreover, our framework includes general random times of
observations, not just doubly stochastic Poisson processes.

We remark that, while being natural, the Frey and Runggaldier model
adopted in this paper does not quite fit into the ``standard''
diffusion or simple point process filtering frameworks (cf.
\cite{KrZa,LSII,Roz1}) and requires more technical tools. In
particular, the general filtering theory for diffusion processes
requires that the diffusion coefficient of the observation process does
not depend on the state process, while in our case the presence of
$\theta_{t}$ in the diffusion coefficient is crucial. The ``standard''
filtering theory for point processes is also not applicable in the
present setting since the observation process $( \tau_{i},S_{\tau_{i}
})  _{i\geq1}$ is a multivariate process (see also Remark \ref{r1}).

It turns out that the resulting filtering equations are simpler than
their counterparts in the case of continuous observations. In the
latter case, the nonlinear filters are described by
infinite-dimensional stochastic differential equations. For example, if
$\theta_{t}$ is a diffusion process, the filtering equations (e.g.,
Kushner filter or Zakai filter) are given by stochastic partial
differential equations (see, e.g., \cite{Roz1}). In contrast, in our
setting, the filtering equation can be reduced to a recursive system of
linked \textit{deterministic} equations of Kolmogorov type. Therefore,
the numerical implementation of the filter is much simpler (see the
follow up paper \cite{CRZ}).

We describe the model in Section~\ref{Sec-2}, state the main results
and examples in Section~\ref{Sec-3}, provide the proofs in
Section~\ref{Section3}, and present more detailed examples in
Section~\ref{Sec-Ex}.

\section{Mathematical model}\label{Sec-2}

\subsection{Risky asset and observation times}

Let us fix a probability space\break $(\Omega,\mathcal{F},\mathsf{P})$
equipped with a filtration $\mathbf{F}=(\mathcal{F}_{t})_{t\geq0}$ that
satisfies the ``usual'' conditions (see, e.g., \cite{LSMar}). All
random processes considered in the paper are assumed to be defined on
$(\Omega,\mathcal{F},\mathsf{P})$ and adapted to $\mathbf{F}$.

It is assumed that there is a risky asset with the price process
$S=(S_{t})_{t\geq0}$ given by the It\^{o} equation
\begin{equation}
dS_{t}=m(\theta_{t})S_{t}\,dt+v(\theta_{t})S_{t}\,dB_{t},\label{2.1b}%
\end{equation}
where $B=(B_{t})_{t\geq0}$ is a standard Brownian motion and $\theta
=(\theta_{t})_{t\geq0}$ is a c\'{a}dl\'{a}g Markov jump-diffusion process in
$\mathbb{R}$ with the generator $\mathcal{L}$. To simplify the discussion, it
is assumed that $m(x)$ and $v(x)$ are measurable bounded functions on
$\mathbb{R}$, the initial condition $S_{0}$ is constant, and $v(x)$ and
$S_{0}$ are positive.

The process $(\theta_{t})_{t\geq0}$ is referred to as the \emph{volatility
process}. It is unobservable, and the only observable quantities are the
values of the log-price process $X_{t}=\log S_{t}$ taken at stopping times
$(\tau_{k})_{k\geq0}$, so that $\tau_{0}=0,\tau_{k}<\tau_{k+1}$ if $\tau
_{k}<\infty,$ and $\tau_{k} \uparrow\infty$ as $k\uparrow\infty.$

In accordance with \eqref{2.1b}, the log-price process is given by
\[
X_{t}=\int_{0}^{t} \bigl(m(\theta_{s})-\tfrac{1}{2}v^{2}(\theta_{s} )
\bigr)\,ds+\int_{0}^{t}v(\theta_{s})\,dB_{s}.
\]
For notational convenience, set $X_{k}:=X_{\tau_{k}}.$ Thus, the observations
are given by the sequence $(\tau_{k},X_{k})_{k\geq0}$.

\begin{remark}[(\textit{Note on the reading sequence})] The reader
interested primarily in applying our results to real data can focus her
attention on Example \ref{examplechain}, which appears to be the most
practical model to work with. That example provides self-contained
formulas for estimating the conditional (filtering) distribution of the
volatility process. We report on the numerical results related to this
example in the follow-up paper~\cite{CRZ}.
\end{remark}

Clearly, the observation process $(\tau_{k},X_{k})_{k\geq0}$ is a
multivariate (marked) point process (see, e.g., \cite{JS,Last}) with
the counting measure
\[
\mu(dt,dy)=\sum_{k\geq1}\mathbf{I}_{\{  \tau_{k}<\infty\}  }
\delta_{{\{\tau_{k},X_{k}\}}}(t,y)\,dt\,dy,
\]
where $\delta_{{\{\tau_{k},X_{k}\}}}$ is the Dirac delta-function on
$\mathbb{R}_{+}\times\mathbb{R}$.

We introduce two filtrations related to $(\tau_{k},X_{k})_{k\ge0}$:
$(\mathcal{G}(n))_{n\ge0}$ and $(\mathcal{G}_{t})_{t\ge0}$, where
\begin{eqnarray*}
\mathcal{G}(n)&:=& \sigma\{(\tau_{k},X_{k})_{k\leq n}\},
\\
\mathcal{G}_{t}&:=& \sigma\bigl(\mu([0,r]\times\Gamma)\dvtx r\leq
t,\Gamma\in \mathcal{B}(\mathbb{R})\bigr),
\end{eqnarray*}
where $\mathcal{B}(\mathbb{R})$ is the Borel $\sigma$-algebra on
$\mathbb{R}$.

It is a standard fact (see Theorem 31 in Chapter III, Section 3 in \cite{JS}) that
\begin{equation}
\mathcal{G}_{\tau_{k}}=\mathcal{G}(k),\qquad k=0,1,\ldots,\label{2.3Gg}%
\end{equation}
and $\{\tau_{k}\}$ is a system of stopping times with respect to
$(\mathcal{G}_{t})_{t\geq0}$.

\begin{remark}\label{r1}
Although $\mathcal{G}_{\tau_{k}}$ contains all the relevant information
carried by the observations obtained up to time $\tau_{k}$, the
filtration $ (\mathcal{G}_{t} )_{t\geq 0}$ provides additional
information between the observation times. To elucidate this point on a
more intuitive level, we note that the length of the time elapsed
between $\tau_{k}$ and $\tau_{k+1}$ carries additional information
about the state of $\theta_{t}$ after $\tau_{k}.$ Specifically, if the
frequency of observations is proportional to the stock's volatility
$v(\theta_{t})$,
$t\in\rrbracket\tau_{k},\tau_{k+1}\llbracket$, the
larger values of $t-\tau_{k}$ might indicate lower values of
$v(\theta_{t})$.
\end{remark}

\subsection{Volatility process}

A more precise description of the volatility process is in order now. Let
$(\mathbb{R},\mathcal{B}(\mathbb{R}))$ and $(\mathbb{R}_{+}\times
\mathbb{R},\mathcal{B}(\mathbb{R}_{+})\otimes\mathcal{B}(\mathbb{R}))$ be
measurable spaces with Borel $\sigma$-algebras. The volatility process
$\theta=(\theta_{t})_{t\geq0}$ is defined by the It\^{o} equation
\begin{equation}
d\theta_{t}=b(t,\theta_{t})\,dt+\sigma(t,\theta_{t})\,dW_{t}+\int_{\mathbb{R}
}u(\theta_{t-},x)(\mu^{\theta}-\nu^{\theta})(dt,dx),\label{2.1}%
\end{equation}
where $W_{t}$ is a standard Wiener process and $\mu^{\theta}=\mu^{\theta
}(dt,dx)$ is a Poisson measure on $(  \mathbb{R}_{+}\times\mathbb{R}%
,\mathcal{B}(  \mathbb{R}_{+})  \otimes\mathcal{B}(
\mathbb{R})  )  $ with the compensator $\nu^{\theta}%
(dt,dx)=K(dx)\,dt$, where $K(dx)$ is a $\sigma$-finite nonnegative
measure on $(  \mathbb{R},\mathcal{B}(  \mathbb{R})  ) $. We assume
that $E\theta_{0}^{2}<\infty$, the functions $b(t,z),\sigma(t,z)$ and
$u(z,x)$ are Lipschitz continuous in $z$ uniformly with respect to
other variables, and
\[
|b(t,z)|^{2}+|\sigma(t,z)|^{2}+\int_{\mathbb{R}}|u(z,x)|^{2}K(dx)\leq
C(1+|z|^{2}).
\]
It is well known that under these assumptions \eqref{2.1} possesses a unique
strong solution adapted to $\mathbf{F}$, and $E\theta_{t}^{2}<\infty$ for any
$t\geq0$.

The generator $\mathcal{L}$ of the volatility process is given by
\begin{eqnarray*}\label{2.10a}
\mathcal{L}f(x) &:=&
b(t,x)f^{\prime}(x)+\tfrac{1}{2}\sigma^{2}(t,x)f^{\prime
\prime}(x)
\\
&&{} +\int_{\mathbb{R}} \bigl(f\bigl(x+u(x,y)\bigr)-f(x)-f^{\prime}
(x)u(x,y) \bigr)K(dy).
\end{eqnarray*}

Before proceeding with the assumptions and main results, we shall
introduce additional notation. Set
\begin{equation}\label{m}
a(s,t)=\int_{s}^{t}\bigl(
m(\theta_{u})-\tfrac{1}{2}v^{2}(\theta_{u})\bigr)\,du
\end{equation}
and
\begin{equation}
\sigma^{2}(s,t)=\int_{s}^{t}v^{2}(\theta_{u})\,du.\label{sigma}%
\end{equation}
For simplicity, it is assumed that $v^{2}(s,t)$ is bounded away from zero. Let
us denote by $\rho_{s,t}(y)$ the density function of the normal distribution
with mean $a(s,t)$ and the variance $\sigma^{2}(s,t)$:
\begin{equation}
\rho_{s,t}(y):=\frac{1}{\sqrt{2\pi}\sigma(s,t)}e^{-{(y-a(s,t))^{2}
}/(2\sigma^{2}(s,t))}.\label{rho}%
\end{equation}
Clearly, $\rho$ is the conditional density of the stock's log-increments
$X_{t}-X_{s}$ given~$\theta$.

Let $\mathcal{F}^{\theta}_{\infty}=(\mathcal{F}_{t}^{\theta})_{t\geq0}$
be the right-continuous filtration generated by $(\theta_{t})_{t\geq0}$
and augmented by $\mathsf{P}$-zero sets from $\mathcal{F}$. Denote by
$G_{k}^{\theta}$ the conditional  distribution of $\tau_{k+1}$ with
respect to~$\mathcal{F}^{\theta}_{\infty}\vee\mathcal{G}(  k) $ (here
and below $\mathcal{F}^{1}\vee\mathcal{F}^{2}$ stands for the
$\sigma$-algebra generated by the $\sigma$-algebras $\mathcal{F}^{1}$
and $\mathcal{F}^{2}$). That is, $G_{k}^{\theta}$ is the distribution of
the time of the next observation, given previous history, and given
$\theta$,
\begin{equation}
G_{k}^{\theta}(  dt)  =\mathsf{P}\bigl(  \tau_{k+1}\in
dt|\mathcal{F}^{\theta}_{\infty}\vee\mathcal{G}(  k) \bigr).
\label{Gktheta}%
\end{equation}
Without loss of generality, we can and will assume that
$G_{k}^{\theta}( dt)  $ is the regular version of the RHS of
(\ref{Gktheta}).

Let $N=(N_{t})_{t\geq0}$ be the counting process with interarrival
times: $\tau_{0}=0$, $(  \tau_{k}-\tau_{k-1})  _{k\geq1},$ that is,
\begin{equation}
N_{t}=\sum_{k\geq1}I(\tau_{k}\leq t).\label{Nt}%
\end{equation}

\subsection{Assumptions}

The following assumptions will be in force throughout the paper:

\renewcommand{\theass}{A.\arabic{ass}}
\setcounter{ass}{-1}
\begin{ass}\label{A.0}
For every $\mathcal{G}$-predictable and a.s. finite stopping time $S$,
\[
\mathsf{P}(N_{S}-N_{S-}\neq0 |\mathcal{G}_{S-})=0\quad\mbox{or}\quad 1.
\]
\end{ass}

\begin{ass}\label{A.1}
The Brownian motion $B$ is independent of $ (\theta,N )$.
\end{ass}

\begin{ass}\label{A.2}
For every $k$, there exists a $\mathcal{G}( k)  $-measurable integrable
random measure $\Phi_{k}$ on $\mathcal{B}(  \mathbb{R}_{+}) $ so that
for almost all $\omega \in\Omega, \Phi_{k}(  [  0,\tau_{k}( \omega) ]
) =0$ and $G_{k}^{\theta}$ is absolutely continuous with respect to
$\Phi_{k}$.
\end{ass}

Denote by $\phi(  \tau_{k},t)  =\phi(  \theta,\tau _{k},t) $ the
Radon--Nikodym derivative of $G_{k}^{\theta}( dt)  $ with respect to
$\Phi_{k}(  dt)  ,$ that is, for almost
every $\omega,$%
\begin{equation}
\phi(  \tau_{k},t)  :=\frac{dG_{k}^{\theta}(  (\tau
_{k},t])  }{d\Phi_{k}(  (\tau_{k},t])  }.\label{fika}%
\end{equation}

Assumption \ref{A.0} is not essential for the derivation of the filter.
However, under this assumption, the structure of the optimal filter is
simpler, and in the practical examples important for this paper, this
assumption holds anyway. In particular, Assumption~\ref{A.0} is
verified if the conditional distribution $G_{k}^{\theta}=\mathsf{P}(
\tau_{k+1}\leq t|\mathcal{F}^{\theta }_{\infty}\vee\mathcal{G}( k)  ) $
is absolutely continuous with respect to the Lebesgue measure or if the
arrival times $\tau_{k}$ are nonrandom (more generally, it holds if the
compensator of the counting process $N_{t}$ is a continuous process).

The following two simple but important examples illustrate
Assumption~\ref{A.2}.

\begin{example}\label{ex:cox copy(1)}
Let ($\tau_{k})_{k\geq0}$ be  the jump times of a doubly stochastic
Poisson process  (Cox process) with the intensity $n(\theta_{t}).$ In
this case,
\[
\mathsf{P}\bigl(\tau_{k+1}\leq t|\mathcal{F}^{\theta}_{\infty}\vee\mathcal{G}%
(  k)  \bigr)= \cases{ 1-e^{-\int_{\tau_{k}}^{t}n(\theta_{s})\,ds},
&\quad  $t\geq\tau_{k}$, \cr\noalign{} 0, &\quad
otherwise.}
\]
Then, one can take $\Phi_{k}(  ds)  =ds$ and $\phi(\tau
_{k},s)=n(\theta_{t})\exp(  -\!\int_{\tau_{k}}^{s}n(  \theta _{u})\,du)
$. If $n(\theta_{t})=n$ is a constant, one could also choose
\[
\Phi_{k}(  ds)  =n\exp\{  n(  \tau_{k}-s) \}\,  ds
\quad\mbox{and}\quad\phi(\tau_{k},s)=1.
\]
\end{example}

\begin{example}
\label{ex:constep0}If the filtering is based on nonrandom observation
times $\tau_{k}$ (e.g., $\tau_{k}=kh$ where $h$ is a fixed time step),
then a natural choice would be $\Phi_{k}(  ds) =\delta_{\{ \tau_{k+1}\}
}( s)\, ds$ and $\phi (\tau_{k},s)=1.$
\end{example}

For practical purposes, $\Phi_{k}(  ds)  $ must be known or easily
computable as soon as the observations $( \tau_{i},X_{i}) _{i\leq
k}$ become available. In contrast, the Radon--Nikodym density $\phi(
\tau_{k})  $ is, in general, a function of the volatility process and
is subject to estimation.

We note that Assumption~\ref{A.2} could be weakened slightly by
replacing $G_{k}^{\theta}$ by a regular version of the conditional
distribution of $\tau_{k+1}$ with respect to
\mbox{$\mathcal{F}_{\tau_{k+1}-}^{\theta}\vee\mathcal{G}(  k)  $}. The latter
assumption would make the proof a little bit more involved and we leave
it to the interested reader.

\section{Main results and introductory examples}\label{Sec-3}

\subsection{Main result}\label{sec-3.1ab}
For a measurable function $f$ on $\mathbb{R}$ with
$E|f(\theta_{t})|<\infty,$ define the conditional expectation estimator
$\pi_{t}(f)$ by
\begin{equation}
\pi_{t}(f):=E \bigl(f(\theta_{t})|\mathcal{G}_{t} \bigr)=\int_{\mathbb{R}}
f(z)\pi_{t}(dz),\label{eq:p0}%
\end{equation}
where $\pi_{t}(dz):=d\mathsf{P}(\theta_{t}\leq z|\mathcal{G}_{t})$ is
the filtering distribution. [Note that we omit the argument
$\theta_{t}$ of $f$ in the estimator $\pi_{t}(f)$.] In the spirit of
the Bayesian approach, it is assumed that the a priori distribution
\[
\pi_{0}(dx)=\mathsf{P}(  \theta_{0}\in dx)
\]
is given.

Let $\sigma\{\theta_{\tau_{k}}\}$ be the $\sigma$-algebra generated by
$\theta_{\tau_{k}}$. For $t>\tau_{k}$, let us define the following
\textit{structure functions}:%
\begin{equation}
\psi_{k}\bigl(f;t,y,\theta_{\tau_{k}}\bigr):=E
\bigl(f(\theta_{t})\rho_{{\tau_{k},t} }(y-X_{k})\phi(\tau_{k},t)
|\sigma \bigl\{\theta_{\tau_{k}} \bigr\}\vee
\mathcal{G}(  k)   \bigr),\label{nado}%
\end{equation}
and its integral with respect to $y$,
\begin{eqnarray}\label{nado1}%
\bar{\psi}_{k}\bigl(f;t,\theta_{\tau_{k}}\bigr)&:=&\int_{\mathbb{R}}\psi_{k}\bigl(
f;t,y,\theta_{\tau_{k}}\bigr)\,dy \nonumber
\\[-8pt]
\\[-8pt]
\nonumber &\hspace*{3pt} =& E \bigl(f(\theta_{t})\phi(\tau_{k} ,t) |\sigma
\bigl\{\theta_{\tau_{k}} \bigr\}\vee\mathcal{G}( k) \bigr),
\end{eqnarray}
where $\rho$ and $\phi$ are given by (\ref{rho}) and (\ref{fika}),
respectively.

If $f\equiv1$, the argument $f$ in $\psi$ and  $\bar{\psi}$ is replaced
by $1.$

Write
\[
\Phi_{k}(\{\tau_{k+1}\}):=\int_{0}^{\infty}I(t=\tau_{k+1})\Phi_{k}(dt),
\]
that is, $\Phi_{k}(\{\tau_{k+1}\})$ is the jump of $\Phi_{k}(dt)$ at
$\tau_{k+1}$.

Finally, for $t\geq\tau_{k}$ and a bounded function $f$, define
\[
\mathcal{M}_{k}(  f;t,\pi_{t})  :=\frac{\pi_{\tau_{k}}( \bar{\psi}_{k}(
f;t)  ) -\pi_{t-}(f)\pi_{\tau_{k}}( \bar{\psi}_{k}( 1;t) )
}{\int_{t}^{\infty}\pi_{\tau_{k} }( \bar{\psi}_{k}( 1;s) ) \Phi_{k}(
ds) }
\]
whenever the numerator is not zero. If the numerator is zero, set
$\mathcal{M}_{k}(  f;t,\pi_{t})  $ to be equal to zero.

The main result of this paper is as follows:

\begin{theorem}
\label{mainthm} Let Assumptions \textup{\ref{A.0}--\ref{A.2}} hold. Then for every
measurable bounded function $f$ in the domain of the generator
$\mathcal{L}$ such that $\int_{0}
^{t}E|\mathcal{L}f(\theta_{s})|\,ds<\infty$ for any $t\geq0,$ the
following system of equations holds:
\begin{longlist}[(2)]
\item[(1)] For every $k=0,1,\ldots,$%
\begin{eqnarray}\label{eq:jump}%
\pi_{\tau_{k+1}}(f)&=&
\frac{\pi_{\tau_{k}}(\psi_{k}(f;t,y))}{\pi_{\tau_{k}}
(\psi_{k}(1;t,y))}\bigg|_{\Bigl\{\matrix{\scriptstyle{\hspace*{-4pt} t=\tau_{k+1}} \cr\noalign{} \scriptstyle{y=X_{k+1}}}\!\! \Bigr\}}
\nonumber
\\[-8pt]
\\[-8pt]
\nonumber &&{} -\mathcal{M}_{k}(  f;t,\pi_{t})  _{ |\{t=\tau_{k+1}
\}}\Phi(  \{ \tau_{k+1}\}  ).
\end{eqnarray}

\item[(2)] For every $k=0,1,\ldots$ and $t\in\,\rrbracket \tau_{k},\tau
_{k+1}\llbracket$,
\begin{equation}
d\pi_{t}(f)=\pi_{t}(\mathcal{L}f)\,dt-\mathcal{M}_{k}(  f;t,\pi_{t}
)  \Phi_{k}(dt).\label{eq:cont}%
\end{equation}
\end{longlist}
\end{theorem}

\subsection{Remarks}\mbox{}

1. Equations (\ref{eq:jump}) and  (\ref{eq:cont}) form a closed system
of equations for the filter~$\pi_{t}(f)$. It is often convenient and
customary (see, e.g., \cite{Roz1,Roz2} and the references
therein) to write a differential equation for a measure-valued process
$H_{t}(  dx)  $ in its variational form, that is, as the related system
of equations for $H_{t}(  f)  $ for all $f$ from a sufficiently rich
class of test functions belonging to the domain of the
operator~$\mathcal{L}.$  In our setting, such a reduction to the variational
form is a necessity, since in some cases the filtering measure
$\pi_{s}(  dx) =\mathsf{P}( \theta_{s}\in dx|\mathcal{G}_{s})  $ may
not belong to the domain of $\mathcal{L}$. However, in the important
examples discussed below, there is no need to resort to the variational
form. The interested reader who is unaccustomed to the variational
approach might benefit from looking first into the examples at the end
of this section and in Section~\ref{Sec-Ex}, where the filtering
equations are written as equations for posterior distributions.

2. The system (\ref{eq:jump}) simplifies considerably if
\begin{equation}
\mathcal{M}_{k}(  f;t,\pi_{t}) _{ |\{t=\tau_{k+1}\}}\Phi(
\{  \tau_{k+1}\}  )  =0\qquad\mbox{for all }k.\label{m0}%
\end{equation}
Obviously, (\ref{m0}) holds if, for all $k, \Phi_{k}(  dt)  $ is
continuous at $t=\tau_{k+1}$, as in the case when $N_{t}$ is a Cox
process. In fact, (\ref{m0}) holds true in many other interesting
cases, even when $\Phi_{k}(  dt)  $ has jumps at all $\tau_{k+1}$, as
in the case of fixed observation intervals (see Example
\ref{ex:constep} below). We note then that the following
\textit{separation principle} holds.

\begin{corollary}
\label{cor:sp} Assume \textup{(\ref{m0})}. Then the filtering at the
observation times $\{\tau_{k}\}_{k\geq1}$ does not require filtering between
them; it is done by the Bayes type recursion\textup{:}
\begin{equation}
\pi_{\tau_{k+1}}(f)=\frac{\pi_{\tau_{k}}(\psi_{k}(f;t,y))}{\pi_{\tau_{k}}
(\psi_{k}(1;t,y))}\bigg|_{ \Bigl\{\matrix{\hspace*{-4pt} \scriptstyle{t=\tau_{k+1}}\cr\noalign{}
\scriptstyle{y=X_{k+1}}}\!\!\Bigr\}}.\label{eq:jump1}%
\end{equation}
\end{corollary}

3. Note that for high-frequency observations, even if condition
(\ref{m0}) is not met, for all practical purposes, it may suffice to
compute the volatility estimates only at the observation times. In that
case, one would only use the relatively simple recursion formula
(\ref{eq:jump}), and disregard equation (\ref{eq:cont}).

4. Clearly, the ``structure functions'' $\psi$ and $\bar {\psi}$ are of
paramount importance for computing the posterior distribution of the
volatility process. We would like to stress that these do not involve
the observations and could be pre-computed ``off-line'' using just the
a priori distribution. Then, ``on-line,'' when the
observations become available, one needs only to plug in the obtained
measurements $(\tau_{k},X_{k}),$ and to compute $\pi_{t}(f)$ by
recursion. This feature is important for developing efficient numerical
algorithms.

 5. Note also that, for almost every $\omega\in\Omega,$ filtering
equation (\ref{eq:cont}) is a \textit{linear deterministic} equation of
Kolmogorov's type, rather than a \textit{nonlinear stochastic} partial
differential equation. The latter is typical of the nonlinear filtering
of diffusion processes. The well-posedness and the regularity
properties of equation (\ref{eq:cont}) are well researched in the
literature on second-order parabolic deterministic integro-differential
equations (see, e.g., \cite{SK,LM,MP}  and the references therein).

\begin{example}[(\textit{Volatility as a Markov chain})]\label{examplechain}
Let us now assume that the counting process is a Cox process with
intensity $n(\theta_{t})$, and take
$\phi(\tau_{k},s)=n(\theta_{t})e^{-\int_{\tau_{k}
}^{s}n(\theta_{u})\,du}$ and $\Phi_{k}(  ds)  =ds.$ Also assume
$\theta=(\theta_{t})_{t\leq T}$ is a homogeneous Markov jump process
taking values in the finite alphabet $\{a_{1},\ldots,a_{M}\}$ with the
intensity matrix $\Lambda=\|\lambda(  a_{i},a_{j})  \|$ and the initial
distribution $p_{q}=\mathsf{P}(\theta_{0}=a_{q})$, $q=1,\ldots,M$.
(This is one of the two models of the state process discussed
in~\cite{FR}.) In this case,
\[
\mathcal{L}f(  \theta_{s})  =\sum_{j}\lambda(  \theta _{s},a_{j}) f(
a_{j}) .
\]
Denote by $\theta_{t }^{j}$ the process $\theta_{t }$ starting from
$a_{j}$, and
\begin{eqnarray*}
p_{ji}(  t)  &:=& \mathsf{P}(  \theta_{t}=a_{i}|\theta_{0}%
=a_{j}),\qquad \pi_{j}(t)=\mathsf{P}(  \theta_{t}=a_{j}%
 |\mathcal{G}_{t})  ,
 \\
r_{ji}(  t,z)  &:=& E \bigl(e^{-\int_{0}^{t}n(\theta_{u}^{j})\,du}
\rho_{_{0,t}}^{j}(z)|\theta_{t}^{j}=a_{i} \bigr),
\end{eqnarray*}
where $\rho_{_{0,t}}^{j}(z)$ is obtained by substituting $\theta_{s}
^{j}$ for $\theta_{s}$ in $\rho_{_{0,t}}(z).$ It follows from Theorem
\ref{mainthm} (for details, see Example \ref{exs1}), with $f(
\theta_{t})  :=I_{\{  \theta_{t}=a_{i}\}  },$ that
\begin{equation}
\hspace*{10mm} \pi_{i}(\tau_{k})=\frac{n(  a_{i})  \sum_{j}r_{ji}(  \tau
_{k}-\tau_{k-1},X_{k}-X_{k-1})  p_{ji}(  \tau_{k}-\tau _{k-1})
\pi_{j}(\tau_{k-1})}{\sum_{i,j}n(  a_{i}) r_{ji}(
\tau_{k}-\tau_{k-1},X_{k}-X_{k-1})  p_{ji}( \tau
_{k}-\tau_{k-1})  \pi_{j}(\tau_{k-1})}.\label{chainTk}%
\end{equation}
This recursion can be easily computed, once one computes (``off-line'')
the values~$r_{ij}$. This example is also treated in more detail in
Section \textup{\ref{Sec-Ex}}.
\end{example}

\section{Proofs}\label{Section3}

In the proof of the main result we want to show that
\[
d\pi_{t}(f) =\pi_{t}(\mathcal{L}f)\,dt +dM_{t},
\]
where $M_{t}$ is a martingale, and then we find a (integral) martingale
representation of $M_{t}$ with respect to the measure $\mu-\nu$, where $\nu$
is a compensator of $\mu$. We first find the compensator.

\subsection{${(\mathcal{G}_{t})}$-compensator of $\mu$}\label{subsec-3}

Denote by $\mathcal{P}(\mathcal{G})$ the predictable
$\sigma$-algebra on $\Omega\times\lbrack0,\infty)$ with respect to
$\mathcal{G}$ and set
\[
\widetilde{\mathcal{P}}(\mathcal{G})=\mathcal{P}(\mathcal{G})\otimes
\mathcal{B}(\mathbb{R}).
\]

A nonnegative random measure $\nu(dt,dy)$ on $\widetilde{\mathcal{P}
}(\mathcal{G})$ is called a \mbox{$\widetilde{\mathcal{P}}(\mathcal{G})$
-}compensa\-tor of $\mu$ if, for any $\widetilde{\mathcal{P}}(\mathcal{G}%
)$-measurable, nonnegative function $\varphi(t,y)=\varphi(\omega,t,y)$:
\begin{longlist}[\hspace*{7mm} (i)]
\item[(i)]
$\displaystyle{ 
\int_{0}^{t}\int_{\mathbb{R}}\varphi(s,y)\nu(ds,dy)
\mbox{ is ${\mathcal{P}}(\mathcal{G})$-measurable},}$

\begin{equation}
\end{equation}

\item[(ii)]
$\displaystyle{E\int_{0}^{\infty}\int_{\mathbb{R}}\varphi
(t,y)\mu(dt,dy)=E\int_{0}^{\infty}\int_{
\mathbb{R}}\varphi(t,y)\nu(dt,dy).}$
\end{longlist}
Let $G_{k}(  ds,dx)  =G_{k}(  \omega,ds,dx)  $ be a regular version of
the conditional distribution of $(  \tau_{k+1} ,X_{k+1})  $ given
$\mathcal{G}(  k)  $ (it is assumed that $G_{k}(  [ 0,\tau_{k}] ,dx)
=0$):
\begin{equation}\label{GGG}
\mathsf{G}_{k}(dt,dy)=d\mathsf{P} \bigl(\tau_{k+1}\le
t,X_{k+1}\le y|\mathcal{G}(k) \bigr).
\end{equation}
Denote $G_{k}(  dt)  =G_{k}(  dt,\mathbb{R})  ,$ that is,
$G_{k}(t)=\mathsf{P}(\tau_{k+1}\leq t | \mathcal{G}(  k) )$ (with
probability one).

By Theorem III.1.33 in \cite{JS}  (see also Proposition 3.4.1 in
\cite{LSMar}),
\begin{equation}
{\nu}(dt,dy)=\sum_{k\geq0}I_{\rrbracket \tau_{k},\tau_{k+1}\rrbracket }
(t)\frac{G_{k}(dt,dy)}{G_{k}([t,\infty),\mathbb{R})}.\label{not:nu}%
\end{equation}
We now derive a representation, suitable for the filtering purposes, of
the \mbox{$\widetilde{\mathcal{P}}(\mathcal{G})$-}compensator $\nu$ in terms
of the structure functions (\ref{nado}), (\ref{nado1}) and the
posterior distribution of~$\theta$.

\begin{lemma}\label{lem-4.1}
The  $\widetilde{\mathcal{P}}(\mathcal{G})$-compensator $\nu$
admits the following version:
\begin{equation}
\nu(dt,dy)=\sum_{k\geq0}I_{\rrbracket \tau_{k},\tau_{k+1}\rrbracket}
(t)\frac{\pi_{\tau_{k}}(\psi_{k}(1;t,y))}{\int_{t}^{\infty}\pi_{\tau_{k}
}(\bar{\psi}_{k}(1;s))\Phi_{k}(ds)}\Phi_{k}(dt)\,dy.\label{eq:nu}%
\end{equation}
\end{lemma}

\begin{pf}
By Assumption~\ref{A.1}, for $t>\tau_{k}$, with probability 1,%
\begin{eqnarray}\label{eq:fteta}%
\nonumber  &&  \mathsf{P} \bigl(\tau_{k+1}\leq t,X_{k+1}\leq
y|\mathcal{F}^{\theta}_{\infty
}\vee\mathcal{G}(  k)   \bigr)
\\
\nonumber &&\qquad  =E \bigl(\mathsf{P} \bigl(\tau_{k+1}\leq t,X_{k+1}\leq y|\mathcal{F}^{\theta
}\vee\mathcal{G}(  k)  \vee\sigma( \tau_{k+1})
 \bigr) |\mathcal{F}^{\theta}_{\infty}\vee\mathcal{G}(  k) \bigr)
\\
&&\qquad  =E\bigl(  I_{(  \tau_{k+1}\leq t)  }\mathsf{P}\bigl( X_{k+1}\leq
y|\mathcal{F}^{\theta}_{\infty}\vee\mathcal{G}( k) \vee\sigma(
\tau_{k+1})  \bigr)\big|\mathcal{F}^{\theta}_{\infty
}\vee\mathcal{G}(  k)  \bigr)
\\
\nonumber &&\qquad  =E\biggl(  I_{(  \tau_{k+1}\leq t)  }\int_{-\infty}^{y}\rho
_{\tau_{k},\tau_{k+1}}(  z-X_{k})\,dz\Big|\mathcal{F}^{\theta}_{\infty
}\vee\mathcal{G}(  k)  \biggr)
\\
\nonumber &&\qquad  =\int_{\tau_{k}}^{t}\int_{-\infty}^{y}\rho_{\tau_{k},s}(z-X_{k}
)\,dz\,G_{k}^{\theta}(  ds)  ,
\end{eqnarray}
where we recall that $G_{k}^{\theta}$ is a regular version of the
conditional  distribution of $\tau_{k+1}$ with respect to
$\mathcal{F}^{\theta}_{\infty }\vee\mathcal{G}(  k)  $. Thus, by Assumption~\ref{A.2},
for $t>\tau_{k}$, with
probability~1,
\begin{eqnarray}\label{eq:sep}
&& \mathsf{P}\bigl(  \tau_{k+1}\leq t,X_{k+1}\leq y|\mathcal{F}^{\theta
}_{\infty}\vee\mathcal{G}(  k)  \bigr)
\nonumber
\\[-6pt]
\\[-6pt]
\nonumber
&&\qquad   =\int_{\tau_{k}}^{t}\int_{-\infty}^{y}\rho_{\tau_{k},s}( z-X_{k})
\phi(\tau_{k},s)\, dz\,\Phi_{k}(  ds) .
\end{eqnarray}
By (\ref{nado}), using notation (\ref{eq:p0}), we see that
\[
E \bigl(E \bigl[\phi(\tau_{k},s)\rho_{\tau_{k},s}(z-X_{k})|\sigma\bigl\{\theta
_{\tau_{k}}\bigr\}\vee\mathcal{G}(  k)   \bigr]|\mathcal{G}( k)
 \bigr)=\pi_{\tau_{k}}\bigl(\psi_{k}(1;s,z)\bigr).
\]
This, together with (\ref{eq:sep}), yields, recalling definition~(\ref{GGG}),
\begin{equation}
G_{k} (ds,dz)=\pi_{\tau_{k}}\bigl(\psi_{k}(1;s,z)\bigr)\Phi_{k}(ds)\,dz.\label{eq:nom}%
\end{equation}
In the same way, for $t>\tau_{k}$, with probability 1,%
\begin{equation}
G_{k} ([t,\infty],\mathbb{R})=\int_{t}^{\infty}\pi_{\tau_{k}}
\bigl(\bar{\psi}_{k}(1;s)\bigr)\Phi_{k}(ds).\label{eq:denom}%
\end{equation}
This completes the proof.
\end{pf}

\begin{remark}\label{re:den0}
If the right-hand side of  (\ref{eq:denom}) is zero, then
\[
\mathsf{P}\bigl(  \tau_{k+1}\geq t|\mathcal{G}(  k)  \bigr) =0.
\]
Hence, $I_{\rrbracket \tau_{k},\tau_{k+1}\rrbracket }(t)=0$ with
probability 1 and, by the $ 0/0=0$ convention, the corresponding term
in \textup{(\ref{eq:nu})} is zero.
\end{remark}

\subsection{Semimartingale representation of the optimal filter}\label{sec-4}
In this section we will prove the following result.

\begin{theorem}
\label{theo-4.1} For any bounded function $f$ from the domain of the operator
$\mathcal{L}$ such that $\int_{0}^{t}E|\mathcal{L}f(\theta_{s})|\,ds<\infty$ for
all $t<\infty$, the It\^o differential of the optimal filter $\pi_{s}(f)$ is given
by equation
\begin{eqnarray}\label{filteq}
\nonumber d\pi_{s}(f)   &=& \pi_{s}(\mathcal{L}f)\,ds
\\[4pt]
&&{}  +\int_{\mathbb{R}} \Biggl(\,\sum_{k\geq0}I_{\rrbracket \tau_{k},\tau
_{k+1}\rrbracket }(s)\frac{\pi_{\tau_{k}}(\psi_{k}(f;s,y))}{\pi_{\tau_{k}
}(\psi_{k}(1;s,y))}-\pi_{s-}(f) \Biggr)
\\[4pt]
\nonumber && \hspace*{10mm} {}\times
(\mu-\nu)(ds,dy).\vadjust{\goodbreak}
\end{eqnarray}
\end{theorem}

\begin{pf}
It suffices to verify the statement for twice continuously differentiable
functions $f$ with $f,f^{\prime}f^{\prime\prime}$ bounded. By It\^{o}'s formula,
\begin{eqnarray*}
f(\theta_{t})&=& f(\theta_{0})+\int_{0}^{t}\mathcal{L}f(\theta_{s})\,ds+\int
_{0}^{t}f^{\prime}(\theta_{s})\sigma(\theta_{s})\,dW_{s}
\\
&&{} +\int_{0}^{t}\int_{R}\bigl(  f\bigl(  \theta_{s-}+u(  \theta _{s-},x) \bigr) -f(
\theta_{s-})  \bigr)  ( \mu^{\theta}-\nu^{\theta}) (ds,dx).
\end{eqnarray*}
Denote
\begin{eqnarray*}
L_{t}&=& \int_{0}^{t}f^{\prime}(\theta_{s})\sigma(\theta_{s})\,dW_{s}
\\
&&{} +\int_{0}^{t}\int_{\mathbb{R}}\bigl(  f\bigl(\theta_{s-}+u(\theta_{s-},x)\bigr)-f( \theta_{s-})
\bigr)  (\mu^{\theta}-\nu^{\theta})(ds,dx).
\end{eqnarray*}
Then, we have
\[
\pi_{t}(f) = E \bigl(f(\theta_{0})|\mathcal{G}_{t} \bigr)
+E \biggl(\int_{0}^{t}\mathcal{L}f(\theta_{s})\,ds \big|\mathcal{G}%
_{t} \biggr)+E (L_{t}|\mathcal{G}_{t} ).
\]
Set
\begin{eqnarray*}
M_{t}  &=&  \bigl\{E \bigl(f(\theta_{0})|\mathcal{G}_{t} \bigr)-\pi_{0}(f) \bigr\}
\\
&&{} + \biggl\{E \biggl(\int_{0}^{t}\mathcal{L}f(\theta_{s})\,ds \Big|\mathcal{G}%
_{t} \biggr)-\int_{0}^{t}\pi_{s} (\mathcal{L}f )\,ds \biggr\}+E (L_{t}%
|\mathcal{G}_{t} ).
\end{eqnarray*}
Obviously, the process $E (f(\theta_{0})|\mathcal{G}_{t} )-\pi_{0}(f)$
is a $\mathcal{G}_{t}$-martingale. Process $L_{t}$ is a $\mathcal{F}_{t}%
$-martingale. Since $\mathcal{G}_{t}\subseteq\mathcal{F}_{t}$, for
$t>t^{\prime},$
\[
E \bigl(E(L_{t}|\mathcal{G}_{t})|\mathcal{G}_{t^{\prime}} \bigr)=E \bigl(E(L_{t}%
|\mathcal{F}_{t^{\prime}})|\mathcal{G}_{t^{\prime}} \bigr)=E(L_{t^{\prime}%
}|\mathcal{G}_{t^{\prime}}).
\]
Consequently, $E(L_{t}|\mathcal{G}_{t})$ is a martingale too. Finally,
$E (\int_{0}^{t}\mathcal{L}f(\theta_{s})\,ds|\mathcal{G}_{t} )-\int
_{0}^{t}\pi_{s} ((\mathcal{L}f) )\,ds$ is also a $\mathcal{G}_{t}%
$-martingale. Indeed, for $t>s>t^{\prime},$ we have $E (\pi_{s}%
 (\mathcal{L}f) |\break \mathcal{G}_{t^{\prime}} )=E (\mathcal{L}%
f(\theta_{s})|\mathcal{G}_{t^{\prime}} )$, which yields
\begin{eqnarray*}
&&  E\biggl[  E \biggl(\int_{0}^{t}\mathcal{L}f(\theta_{s})\,ds \Big|\mathcal{G}%
_{t} \biggr)-\int_{0}^{t}\pi_{s}(\mathcal{L}f)\,ds \Big|\mathcal{G}_{t^{\prime}%
}\biggr]
\\
&&\qquad  =E \biggl(\int_{0}^{t^{\prime}}\mathcal{L}f(\theta_{s}%
)\,ds \Big|\mathcal{G}_{t^{\prime}} \biggr)-\int_{0}^{t^{\prime}}\pi_{s}%
(\mathcal{L}f)\,ds.
\end{eqnarray*}
Thus, $M_{t}$ is a $\mathcal{G}_{t}$-martingale. In particular, this means
that $\pi_{t}(f)$ is a \mbox{$\mathcal{G}$-}semi\-martin\-gale with paths in the
Skorokhod space $\mathbb{D}_{[0,\infty)}(\mathbb{R})$, so that $\pi_{t}(f)$ is
a right continuous process with limits from the left. By the martingale
representation theorem (see, e.g., Theorem 1 and Problem 1.c in Chapter~4, Section~8 in
\cite{LSMar}),
\[
M_{t}=\int_{0}^{t}\int_{\mathbb{R}}H(s,y)(\mu-\nu)(ds,dy). %
\]
It is a standard fact that $\mathsf{P}(N_{S}-N_{S-}\neq0|\mathcal{G}%
_{S-})={\nu}(\{  S\}  ,\mathbb{R}_{+}).$ Hence, due to Assumption~\ref{A.0},
by Theorem 4.10.1 from \cite{LSMar} [see formulae (10.6)~and~(10.15)],
\begin{equation}
H(t,y)=\mathsf{M}_{\mu}^{\mathsf{P}} \bigl(\triangle M|\widetilde{\mathcal{P}%
}(\mathcal{G}) \bigr)(t,y),\label{4.7}%
\end{equation}
where $\triangle M_{t}=M_{t}-M_{t-}$ and the conditional expectation
$\mathsf{M}_{\mu}^{\mathsf{P}} (g|\widetilde{\mathcal{P}}(\mathcal{G}%
) )$ is defined by the following relation (see, e.g., \cite{LSMar},
Chapter~2, Section~2 and Chapter~10, Section~1): for any $\widetilde{\mathcal{P}}(\mathcal{G}%
)$-measurable bounded and compactly supported function $\varphi(t,y),$
\begin{eqnarray*}\label{4.10}
&& E\int_{0}^{\infty}\int_{\mathbb{R}}\varphi(t,y)g_{t}\mu(dt,dy)
\\
&&\qquad =E\int_{0}^{\infty}\int_{\mathbb{R}}\varphi(t,y)\mathsf{M}_{\mu}^{\mathsf{P}%
} \bigl(g |\widetilde{\mathcal{P}}(\mathcal{G}) \bigr)(t,y)\nu(dt,dy).
\end{eqnarray*}
By Lemma 4.10.2, \cite{LSMar},
\begin{equation}
\mathsf{M}_{\mu}^{P} \bigl(\pi_{t}(f) |\widetilde{\mathcal{P}}(\mathcal{G}%
) \bigr)(t,y)=\mathsf{M}_{\mu}^{P} \bigl(f |\widetilde{\mathcal{P}}%
(\mathcal{G}) \bigr)(t,y).\label{4.10a}%
\end{equation}
Since $\pi_{t-}(f)$ is
$\widetilde{\mathcal{P}}(\mathcal{G})$-measurable [which implies
$\mathsf{M}_{\mu}^{\mathsf{P}}(\pi_{-}(f)|\widetilde
{\mathcal{P}}(\mathcal{G}))(t,y)=\pi_{t-}(f)$], by (\ref{4.10a}),
\begin{eqnarray}\label{4.9}
\nonumber && \mathsf{M}_{\mu}^{\mathsf{P}} \bigl(\triangle M |\widetilde{\mathcal{P}%
}(\mathcal{G}) \bigr)(t,y)
\\
&&\qquad =\mathsf{M}_{\mu}^{\mathsf{P}} \bigl(\pi_{t}(f)-\pi_{t-}(f) |\widetilde
{\mathcal{P}}(\mathcal{G}) \bigr)(t,y)
\\
\nonumber &&\qquad  =\mathsf{M}_{\mu}^{\mathsf{P}} \bigl(f |\widetilde{\mathcal{P}%
}(\mathcal{G}) \bigr)(t,y)-\pi_{t-}(f).
\end{eqnarray}
To complete the proof, one needs to show that
\begin{equation}
\mathsf{M}_{\mu}^{\mathsf{P}} \bigl(f(  \theta_{\bolds{.}})  |\widetilde
{\mathcal{P}}(\mathcal{G}) \bigr)(s,y)=\sum_{k\geq0}I_{\rrbracket \tau
_{k},\tau_{k+1}\rrbracket }(s)\frac{\pi_{\tau_{k}}(\psi_{k}(f;s,y))}%
{\pi_{\tau_{k}}(\psi_{k}(1;s,y))}.\label{gg}%
\end{equation}

To prove (\ref{gg}), it suffices to demonstrate that, for any $\widetilde
{\mathcal{P}}(\mathcal{G})$-measurable bounded and compactly supported
function $\varphi(t,y),$
\begin{eqnarray}\label{show}
\nonumber && E\sum_{k\geq0}\int_{(\tau_{k},\tau_{k+1}]\cap(\tau_{k},\infty)}\int
_{\mathbb{R}}\varphi(t,y)\frac{\pi_{\tau_{k}}(\psi_{k}(f;t,y))}{\pi_{\tau_{k}
}(\psi_{k}(1;t,y))}\nu(dt,dy)
\nonumber
\\[-8pt]
\\[-8pt]
\nonumber
&&\qquad =E\int_{0}^{\infty}\int_{\mathbb{R}} \varphi(t,y)f(\theta_{t})\mu
(dt,dy).
\end{eqnarray}

By monotone class arguments, we can assume that $\varphi(  t,x) =v( t)
g(  x)  $, where $v(  t)  $ is a $\mathcal{P}(\mathcal{G})$-measurable
process and $g( x)  $ is a continuous function on $\mathbb{R}$. By
Lemma~III.1.39 in~\cite{JS}, since ${v}(t)$ is
$\mathcal{P}(\mathcal{G})$-measurable, it must be of the form
\begin{equation}
v(t)  =v_{0}+\sum_{k\geq1}^{\infty}v_{k}(  t)
I_{\rrbracket \tau_{k},\tau_{k+1}\rrbracket }(  t),\label{4.11b}%
\end{equation}
where $v_{0}$ is a constant and $v_{k}(  t)  $ are $\mathcal{G} ( k)
\otimes\mathcal{B}(  \mathbb{R}_{+}) $-measurable functions.

Owing to (\ref{4.11b}) and Lemma \ref{lem-4.1}, in order to prove
(\ref{show}), it suffices to verify the equality
\begin{eqnarray}\label{enough}%
&& E\biggl[  \int_{(\tau_{k},\tau_{k+1}]\cap(\tau_{k},\infty)}\int_{\mathbb{R}
}g(y)v_{k}(t)\frac{\pi_{\tau_{k}}(\psi_{k}(f;t,y))}{\pi_{\tau_{k}}(\psi
_{k}(1;t,y))}\Phi_{k}(  dt)\,dy\biggr]
\nonumber
\\[-8pt]
\\[-8pt]
\nonumber
&&\qquad =E\bigl[  v_{k}(\tau_{k+1})g(X_{k+1})f\bigl(\theta_{\tau_{k+1}}\bigr)\mathbf{1}_{\{
\tau_{k+1}<\infty\}  }\bigr] .
\end{eqnarray}

The next step follows the ideas of Theorem III.1.33 in \cite{JS}. We have
\begin{eqnarray*}
&&  E\bigl[  v_{k}(\tau_{k+1})g(X_{k+1})f\bigl(\theta_{\tau_{k+1}}\bigr)\mathbf{1}_{\{
\tau_{k+1}<\infty\}  }\bigr]
\\
&&\qquad  =E\bigl[  E\bigl(  v_{k}(\tau_{k+1})g(X_{k+1})f\bigl(\theta_{\tau_{k+1} }\bigr)\mathbf{1}_{\{
\tau_{k+1}<\infty\}  }|\mathcal{G}(  k)
\vee\mathcal{F}^{\theta}_{\infty}\bigr)  \bigr]
\\
&&\qquad  =E\biggl(  \int_{(\tau_{k},\infty)}\int_{\mathbb{R}}v_{k}(s)g(y)E[
f(\theta_{s})G_{k}^{\theta}(  ds,dy)  |\mathcal{G}( k) ] \biggr),
\end{eqnarray*}
where, as before, $G_{k}^{\theta} (ds,dy)$ is the regular version of
the conditional  distribution of $(  \tau_{k+1},X_{k+1})  $ with
respect to $\mathcal{F}^{\theta}_{\infty}\vee\mathcal{G}( k)  .$

By Fubini's theorem, and recalling notation \textup{(\ref{GGG})},
\begin{eqnarray}\label{eq:fub}
\nonumber  &&  E\biggl(  \int_{(\tau_{k},\infty)}\int_{\mathbb{R}}v_{k}(s)g(y)E[
f(\theta_{s})G_{k}^{\theta}(  ds,dy)  |\mathcal{G}( k)
]  \biggr)
\\
&&\qquad  =E\biggl(  \int_{(\tau_{k},\infty)}\int_{\mathbb{R}}v_{k}(s)g(y)\frac {E[
f(\theta_{s})G_{k}^{\theta}(  ds,dy) |\mathcal{G}( k) ] }{G_{k}( [
s,\infty] ;\mathbb{R})}
\nonumber
\\[-8pt]
\\[-8pt]
\nonumber
&&\hspace*{58.8mm} {}\times \int_{[s,\infty]}G_{k}(  du,\mathbb{R})  \biggr)
\\
\nonumber &&\qquad  =E\biggl(  \int_{\tau_{k}}^{\tau_{k+1}}\int_{\mathbb{R}}v_{k}(s)g(y)\frac
{E[  f(\theta_{s})G_{k}^{\theta}(  ds,dy) |\mathcal{G}( k) ] }{G_{k}( [
s,\infty] ;\mathbb{R}) }\biggr).
\end{eqnarray}
By (\ref{eq:sep}),
\begin{equation}
G_{k}^{\theta}(  ds,dy)  =\rho_{{\tau_{k},s}}(z-X_{k})\phi(\tau
_{k},s)\Phi_{k}(  ds)\,dy.\label{2.1bb}%
\end{equation}

Hence, for $s>\tau_{k}$,
\begin{eqnarray*}
&&  E[  f(\theta_{s})G_{k}^{\theta}(  ds,dy)  |\mathcal{G}
(  k)  ]
\\
&&\qquad =E \bigl(E\bigl(  f(\theta_{s})\rho_{{\tau_{k},s}}(y-X_{k})\phi(\tau
_{k},s)|\sigma\bigl\{  \theta_{\tau_{k}}\bigr\}  \vee\mathcal{G}(
k)  \bigr)   |\mathcal{G}(k) \bigr)\Phi_{k}(  ds)\,dy.
\\
&&\qquad  =\pi_{\tau_{k}}\bigl(\psi_{k}(  f;s,y)  \bigr)\,dy\,\Phi_{k}(
ds).
\end{eqnarray*}

This, together with \textup{(\ref{eq:denom})}, yields%
\begin{eqnarray*}
&&  E\biggl(  \int_{\tau_{k}}^{\tau_{k+1}}\int_{\mathbb{R}}v_{k}(s)g(y)\frac
{E[  f(\theta_{s})G_{k}^{\theta}(  ds,dy) |\mathcal{G}( k) ] }{G_{k}( [
s,\infty] ;\mathbb{R})
}\biggr)
\\
&&\qquad  =E\biggl(  \int_{\tau_{k}}^{\tau_{k+1}}\int_{\mathbb{R}}v_{k}(s)g(y)\frac
{\pi_{\tau_{k}}(\psi_{k}(  f;s,y)  )\,dy}{\int_{s}^{\infty}\pi
_{\tau_{k}}(  \bar{\psi}(  1;t)  ) \Phi_{k}( dt) }\Phi_{k}( ds) \biggr) ,
\end{eqnarray*}
so that \textup{(\ref{enough})} is satisfied, and the proof follows.
\end{pf}

\subsection[Proof of Theorem 3.1]{Proof of Theorem \textup{\protect\ref{mainthm}}}\label{sec-6}
In this section we show that Theorem \ref{mainthm} follows from
Lemma \ref{lem-4.1} and Theorem \ref{theo-4.1}.

\begin{pf}
First, we note that the stochastic integral in the RHS of (\ref{filteq}) can
be written as the difference of the integrals with respect to $\mu$ and $\nu.$
Indeed, since $f$~is bounded, this follows from \cite{JS}, Proposition II.1.28.

By applying Lemma \ref{lem-4.1} and integrating over $y$, one gets that, for
$t\in\,\rrbracket \tau_{k},\tau_{k+1}\rrbracket ,$%
\begin{eqnarray*}
&& \int_{\mathbb{R\times}(\tau_{k},t]} \biggl(\frac{\pi_{\tau_{k}}(\psi
_{k}(f;s,y))}{\pi_{\tau_{k}}(\psi_{k}(1;s,y))}-\pi_{s-}(f) \biggr)\nu(ds,dy)
\\
&&\qquad  =\int_{(\tau_{k},t]}\frac{\pi_{\tau_{k}}( \bar{\psi}_{k}(f;s))
-\pi_{s-}(f)\pi_{\tau_{k}}( \bar{\psi}_{k}(1;s))  }{\int
_{s}^{\infty}\pi_{\tau_{_{k}}}( \bar{\psi}_{k}(1;u))  \Phi _{k}( du)
}\Phi_{k}( ds).
\end{eqnarray*}

This equation verifies that \eqref{eq:cont} follows from the semimartingale
representation~\eqref{filteq}, for $t$ between the consecutive observation times.

For the jump part \eqref{eq:jump}, we note that
\[
\int_{0}^{t}\int_{\mathbb{R}}\pi_{s-}(  f)  \mu( ds,dy)
=\sum_{\tau_{k+1}\leq t}\pi_{(  \tau_{k+1}) -}( f)
\]
and
\[
\int_{0}^{t}\int_{\mathbb{R}}\frac{\pi_{\tau_{k}}(\psi_{k}(f;s,y))}{\pi
_{\tau_{k}}(\psi_{k}(1;s,y))}\mu(  ds,dy)  =\sum_{\tau_{k+1}\leq
t}\frac{\pi_{\tau_{k}}(\psi_{k}(f;s,y))}{\pi_{\tau_{k}}(\psi_{k}(1;s,y))}
\bigg|_{  \Bigl\{\matrix{\scriptstyle{\hspace*{-4pt} s=\tau_{k+1}}\cr\noalign{}
\scriptstyle{y=X_{k+1}}}\!\!\Bigr\}}.
\]
Now, (\ref{filteq})   can be rewritten as follows:
\begin{eqnarray}\label{closeq1}
\nonumber \pi_{t}(  f)   &=& \pi_{0}(  f)  +\int_{0}^{t}\pi
_{s}(  \mathcal{L}f)\,ds
\\
&&{} +\sum_{\tau_{k+1}\leq t}\biggl(  \frac{\pi_{\tau_{k}}(\psi_{k}(f;s,y))}
{\pi_{\tau_{k}}(\psi_{k}(1;s,y))} \bigg|_{  \Bigl\{
\matrix{\scriptstyle{\hspace*{-4pt} s=\tau_{k+1}}\cr\noalign{} \scriptstyle{y=X_{k+1}}}\!\! \Bigr\}} -\pi_{(  \tau_{k+1})  -}(  f)
\biggr)
\\
\nonumber &&{}  -\sum_{k\geq0}\int_{(\tau_{k},t\wedge\tau_{k+1}]}\mathcal{M}_{k}(
f;s,\pi_{s})  \Phi_{k}(  ds).
\end{eqnarray}
Suppose $t\in\,\rrbracket \tau_{k},\tau_{k+1}\llbracket.$ Then,
\[
\pi_{t}(  f)    = \pi_{\tau_{k}}(  f)
 +\int_{\tau_{k}}^{t}\pi_{s}(  \mathcal{L}f)\,ds-\int_{\tau_{k}
}^{t}\mathcal{M}_{k}(  f;s,\pi_{s})  \Phi_{k}( ds)  .
\]
It follows that
\begin{eqnarray*}
&& \pi_{(  \tau_{k+1})  -}(  f)
\\
&&\qquad  =\pi_{\tau_{k}}(  f)  +\int_{\tau_{k}}^{\tau_{k+1}}\pi _{s}(
\mathcal{L}f)\,ds-\int_{\tau_{k}}^{(  \tau _{k+1}) -}\mathcal{M}_{k}(
f;s,\pi_{s}) \Phi_{k}( ds)  .
\end{eqnarray*}
Therefore, from \textup{(\ref{closeq1})},
\[
\pi_{\tau_{k+1}}(  f)  =\frac{\pi_{\tau_{k}}(\psi_{k}(f;s,y))}
{\pi_{\tau_{k}}(\psi_{k}(1;s,y))} \bigg|_{  \Bigl\{\matrix{\scriptstyle{\hspace*{-4pt}
s=\tau_{k+1}}\cr\noalign{}
\scriptstyle{y=X_{k+1}}}\!\! \Bigr\}} -\mathcal{M}_{k}(  f;t,\pi_{t})
_{ |\{t=\tau_{k+1}\}} \Phi(  \{  \tau_{k+1}\}  )  .
\]
This completes the proof.
\end{pf}

\section{Examples}\label{Sec-Ex}

In this section we consider some important special cases of Theorem~\ref{mainthm}.

\begin{example}[(\textit{Markov chain volatility and Cox process arrivals})]
\label{exs1} Recall the setting of Example~\ref{examplechain} and its
notation $r_{ij}$, $\pi_{j}(t)$ and $\theta^{j}$. It follows from
Example~\ref{ex:cox copy(1)} that in this case $\Phi_{k}(  \{
\tau_{k+1}\}  )  =0$ for all $k.$ Hence, the second term in the RHS
of~\textup{(\ref{eq:jump})} is zero. By (\ref{nado}), for
$f(\theta _{t})=\mathbf{1}_{\{\theta_{t}=a_{i}\}}$ and $t>\tau_{k}$,
\[
\psi_{k}\bigl(f;t,y,\theta_{\tau_{k}}\bigr)=n(  a_{i})  \bigl[ E \bigl(I_{\{
\theta_{t}=a_{i}\} }e^{-\int_{s}^{t}n(\theta_{u}
)\,du}\rho_{_{s,t}}(y-x)|\theta_{s} \bigr)\bigr] _{ \Bigl\{\matrix{\scriptstyle{\hspace*{-4pt} s=\tau
_{k}}\cr\noalign{}\scriptstyle{x=X_{k}}}\!\! \Bigr\}}.
\]
Thus, owing to the homogeneity of $\theta_{t},$ for $t >\tau_{k}$,
\begin{eqnarray*}
&&  \pi_{\tau_{k}}\bigl(  \psi_{k}(f;t,y)\bigr)
\\
&&\qquad =\sum_{j}n(a_{i})E \bigl(I_{\{  \theta_{t}=a_{i}\}  }e^{-\int
_{s}^{t}n(\theta_{u})\,du}\rho_{s,t}(y-x) |\theta_{s}=a_{j}  \bigr)_{
\Bigl\{\matrix{\scriptstyle{\hspace*{-4pt} s=\tau_{k}}\cr\noalign{}
\scriptstyle{x=X_{k}}}\!\! \Bigr\}}\pi_{j}(\tau_{k})
\\
&&\qquad =\sum_{j}n(a_{i})E \bigl(I_{\{  \theta_{t-s}^{j}=a_{i}\}
}e^{-\int_{0}^{t-s}n(\theta_{u})\,du}\rho_{0,t-s}^{j} (y-x) \bigr)_{
\Bigl\{\matrix{\scriptstyle{\hspace*{-4pt} s=\tau_{k}}\cr\noalign{}
\scriptstyle{x=X_{k}}}\!\! \Bigr\}}\pi_{j}(\tau_{k})
\\
&&\qquad =\sum_{j}n(a_{i})E \bigl[I_{\{  \theta_{t-s}^{j}=a_{i}\} }E
\bigl(e^{-\int_{0}^{t-s}n(\theta_{u})\,du}\rho_{0,t-s}^{j}(y-x) |\theta
_{t-s}^{j} \bigr) \bigr]_{\Bigl\{\matrix{\scriptstyle{\hspace*{-4pt} s=\tau_{k}}\cr\noalign{}
\scriptstyle{x=X_{k}}}\!\! \Bigr\}}
\pi_{j}(\tau_{k})
\\
&&\qquad =\sum_{j}n(a_{i})r_{ji}(  t-\tau_{k},y-X_{k})  p_{ji}( t-\tau_{k})
\pi_{j}(\tau_{k}).
\end{eqnarray*}
Similar formula holds for the denominator of the first term of the RHS
of the equation. Now equation \eqref{chainTk}  follows from~\eqref{eq:jump}.

Mimicking the previous calculations and using the notation
\[
\bar{r}_{ji}(  t)  :=E \bigl(e^{-\int_{0}^{t}n(\theta_{u}^{j}
)\,du}|\theta_{t}^{j}=a_{i} \bigr),
\]
it is readily checked that, for $t>\tau_{k}$,
\[
\pi_{\tau_{k}} \bigl(\bar{\psi}_{k}\bigl(\mathbf{1}_{\{\theta_{t}=a_{i} \}};t\bigr)
\bigr)=n(  a_{i})  \sum_{j}\pi_{j}(\tau_{k})\bar{r} _{ji}(t-\tau_{k})p_{ji}(
t-\tau_{k})
\]
and
\[
\pi_{\tau_{k}} \bigl(\bar{\psi}_{k}(1,t) \bigr)=\sum_{i,j}\pi_{j}(\tau _{k})n(
a_{i})  \bar{r}_{ji}(t-\tau_{k})p_{ji}(  t-\tau _{k})  ,
\]
which are needed in computing \eqref{eq:cont}. It is easily verified that, in
the setting of this example, equation \eqref{eq:cont} reduces to the
following:
\begin{eqnarray}\label{eq:ko2}
d\pi_{i}(t)&=& \sum_{j}\lambda(  a_{j},a_{i})  \pi_{j}(t)\,dt
\nonumber
\\[-8pt]
\\[-8pt]
\nonumber
&&{} +\bar {D}(
\tau_{k},t)  \pi_{i}(t)\,dt+D_{i}( \tau_{k},t)\,dt,%
\end{eqnarray}
where
\begin{eqnarray*}
D_{i}(  \tau_{k},t)   &=& -\frac{n(  a_{i})  \sum
_{j}\bar{r}_{ji}(t-\tau_{k})p_{ji}(  t-\tau_{k}) \pi_{j}(\tau
_{k})}{\int_{t}^{\infty}\sum_{i,j}n(  a_{i}) \bar{r}_{ji}
(s-\tau_{k})p_{ji}(  s-\tau_{k})  \pi_{j}(\tau_{k})\,ds},
\\
\bar{D}(  \tau_{k},t)   &=& \frac{\sum_{l,j}n(  a_{l})
\bar{r}_{jl}(t-\tau_{k})p_{jl}(  t-\tau_{k}) \pi_{j}(\tau_{k}
)}{\int_{t}^{\infty}\sum_{i,j}n(  a_{i}) \bar{r}_{ji}(s-\tau
_{k})p_{ji}(  s-\tau_{k})  \pi_{j}(\tau_{k})\,ds}.
\end{eqnarray*}
Note that equation (\ref{eq:ko2}) is considered for a fixed $\omega$
and $t>\tau_{k}(  \omega)  .$ Therefore,
$\tau_{k}$~and~$\pi_{\bolds{\cdot}}(\tau_{k})$ should be viewed as known quantities.
\end{example}

\begin{example}[(\textit{Poisson arrivals})]\label{Poisson}
Let $\theta$ be still the
same as in Example \ref{exs1}. Suppose that the interarrival times
between the observations are exponential with constant intensity
$n(\theta)\equiv\lambda$. In other words, $N_{t}$ is a Poisson process
with constant parameter $\lambda.$ In this case, the volatility process
$\theta$ is independent of $N_{t}.$ { Then, on the interval
$\tau_{k}<t<\tau_{k+1} $, equation }$( \ref{eq:ko2})  $ reduces to{
\begin{eqnarray}\label{eq:po}
d\pi_{i}(t)  &=& \sum_{j}\lambda(  a_{j},a_{i})  \pi_{j}
(t)\,dt
\nonumber
\\[-8pt]
\\[-8pt]
\nonumber
&&{} -\lambda\Biggl(\sum_{j}p_{ji}(  t-\tau_{k})  \pi_{j}(\tau_{k})-\pi
_{i}(t)\Biggr)\,dt.
\end{eqnarray}
On the other hand, owing to the independence of }$N$ and $\theta,$ {it
is readily checked that on the interval $\tau_{k}<t<\tau_{k+1},$}%
\[
\pi_{i}(t)=\sum_{j}p_{ji}(  t-\tau_{k})  \pi_{j}(\tau_{k}).
\]
Therefore, the filtering equation  (\ref{eq:po}) is simply the forward
Kolmogorov equation for $\theta.$
\end{example}

A similar effect appears also in the following example.

\begin{example}[(\textit{Fixed observation intervals})]\label{ex:constep}
Assume for simplicity that the Markov process $\theta_{t}$ is homogeneous. Also assume
that $\tau_{k}=kh,$ where $h$ is a fixed\vadjust{\goodbreak} time step. Notice that
\begin{equation}
\mathcal{G}_{t}=\mathcal{G}(  k)\qquad \mbox{for any $t\in
\llbracket\tau_{k},\tau_{k+1}\llbracket$}.\label{bom}%
\end{equation}
Denote by $P(t,x,dy)$ the transition probability kernel of the process
$\theta_{t}$, given that $\theta_{0}=x$, and let $T_{t}$ denote the associated
transition operator.

In accordance with Example \ref{ex:constep0}, one can take
\[
\phi
(\tau_{k},t)\equiv1\quad\mbox{and}\quad \Phi_{k}(dt)=\delta_{\{\tau_{k+1}\}}(t)\,dt.
\]
Thus, we
get
\begin{eqnarray}
\psi_{k}\bigl(  f;t,y,\theta_{\tau_{k}}\bigr)   &=& E\bigl[  f(\theta
_{t})\rho_{{\tau_{k},t}}(y-X_{k}) |\sigma\bigl\{  \theta_{\tau_{k}
}\bigr\}  \vee\mathcal{G}(  k)  \bigr]  ,
\\
\bar{\psi}_{k}\bigl(  f;t,\theta_{\tau_{k}}\bigr)   &=& T_{t-\tau_{k} }f\bigl(
\theta_{\tau_{k}}\bigr)  :=\int f(y)\mathsf{P}\bigl( t-\tau_{k}
,\theta_{\tau_{k}},dy\bigr).
\end{eqnarray}
Since $\Phi_{k}(dt)=0$ on $\llbracket\tau_{k},\tau
_{k+1}\llbracket$, \eqref{eq:cont} is reduced to the forward Kolmogorov
equation
\[
\frac{\partial_{t}}{\partial t}\pi_{t}(f)=\pi_{t}(\mathcal{L}f),
\]
subject to the initial condition $\pi_{\tau_{k}}(f).$ The unique solution of
this equation is given by $\pi_{t}(f)=\pi_{\tau_{k}}(T_{t-\tau_{k}}f)$,
$t<\tau_{k+1}$. Hence,
\begin{equation}
\pi_{\tau_{k+1}-}(f)=\pi_{\tau_{k}}(  T_{h}f). \label{minus}%
\end{equation}
Since $\phi(\tau_{k},t)\equiv1,$ the
denominator of $\mathcal{M}_{k}$ is equal to
1 when $t=\tau_{k+1}$. This
together with the formula $\Phi(\{\tau_{k+1}\})=1$ yields
\begin{equation}
\mathcal{M}_k(f;t,\pi_t)|
_{{t=\tau}_{{k+1}}}{\Phi(\{\tau_{k+1}\})=\pi
}_{{\tau}_{{k}}}(  T_{h}f )
-\pi_{\tau_{k+1}-}(f).\label{eq:dopjump}%
\end{equation}
Owing to (\ref{eq:dopjump}), we get $\mathcal{M}_{k}(f;t,\pi
_{t})|_{t=\tau_{k+1}} \Phi
(\{\tau_{k+1}\})=0.$

This yields the following recursion formula:
\begin{eqnarray*}
\pi_{\tau_{k+1}}(  f)   &=& \frac{\pi_{\tau_{k}}(  \psi _{k}( f;t,y) )
}{\pi_{\tau_{k}}( \psi_{k}(1;t,y)  )  }\bigg|_{t=\tau_{k+1},y=X_{\tau_{k+1}}}
\\
&=& \frac{\int_{\mathbb{R}}E(  f(\theta_{t-\tau_{k}})\rho_{0{,t-}\tau
_{k}}(y-z)|\theta_{0}=z)  \pi_{\tau_{k}}(  dz)  } {\int_{\mathbb{R}}E(
\rho_{0{,t-}\tau_{k}}(y-z)|\theta_{0}=z) \pi_{\tau_{k}}( dz)
}\bigg|_{t=\tau_{k+1}, y=X_{\tau_{k+1}}}.
\end{eqnarray*}
\end{example}

\section*{Acknowledgments}
We are grateful to the anonymous Associate Editor
and the referee for their constructive suggestions, especially regarding a
simplified presentation of the results. We are very much indebted to Remigijus
Mikulevicius for many important suggestions, and to Ilya Zaliapin, whose
numerical experiments helped to discover an error in a preprint version of the
paper.\vadjust{\goodbreak}

\printaddresses


\begin{thebibliography}{99}
\bibitem{CRZ}
\textsc{Cvitani\'c, J., Rozovskii, B.} and \textsc{Zaliapin, Il.} (2006). Numerical
estimation of volatility values from discretely observed diffusion data.
\textit{Journal of Computational Finance}. To appear.


\bibitem{EHJ}
\textsc{Elliott, R. J., Hunter, W. C.} and \textsc{Jamieson, B. M.} (1998). Drift and
volatility estimation in discrete time. \emph{J.~Econom. Dynam.
Control} \textbf{22} 209--218.
\MR{1488361}\vadjust{\goodbreak}

\bibitem{FPS}
\textsc{Fouque, J.-P., Papanicolaou, G.} and \textsc{Sircar, R.} (2000).
\emph{Derivatives in Financial Markets with Stochastic Volatility}.
Cambridge Univ. Press.
\MR{1768877}

\bibitem{F}
\textsc{Frey, R.} (1997). Derivative asset analysis in models with level-dependent
and stochastic volatility. \emph{CWI Quarterly} \textbf{10} 1--34.
\MR{1472800}

\bibitem{FR}
\textsc{Frey, R.} and \textsc{Runggaldier, W.} (2001). A nonlinear filtering approach to
volatility estimation with a view towards high frequency data.
\emph{Internat. J. Theoret. Appl. Finance} \textbf{4} 199--210.
\MR{1831267}

\bibitem{GT}
\textsc{Gallant, A. R.} and \textsc{Tauchen, G.} (1998). Reprojecting partially observed
systems with application to interest rate diffusions. \emph{J.~Amer.
Statist. Assoc.} \textbf{93} 10--24.


\bibitem{Gou}
\textsc{Gourieroux, C.} (1997). \emph{ARCH Models and Financial Applications}.
Springer, New York.
\MR{1439744}



\bibitem{JS}
\textsc{Jacod, J.} and \textsc{Shiryaev, A. N.} (1987). \emph{Limit Theorems for Stochastic
Processes.} Springer, Berlin.
\MR{0959133}

\bibitem{KalStr}
\textsc{Kallianpur, G.} and \textsc{Striebel, C.} (1969). Stochastic differential
equations occurring in the estimation of continuous parameter stochastic
processes. \textit{Teor. Veroyatnost. i Primenen.} \textbf{14} 597--622.
\MR{0264780}

\bibitem{KX}
\textsc{Kallianpur, G.} and \textsc{Xiong, J.} (2001). Asset pricing with stochastic
volatility. \emph{Appl. Math. Optim.} \textbf{43} 47--62.
\MR{1804394}

\bibitem{SK}
\textsc{Krein, S. G.} (1982). \emph{Linear Equations in Banach Spaces}.
Birkh\"auser, Boston.
\MR{0684836}

\bibitem{KrZa}
\textsc{Krylov, N. V.} and \textsc{Zatezalo, A.} (2000). Filtering of finite-state
time-non homogeneous Markov processes, a direct approach. \emph{Appl.
Math. Optim.} \textbf{42} 229--258.
\MR{1795610}

\bibitem{Last}
\textsc{Last, G.} and \textsc{Brandt, A.} (1995). \textit{Marked Point Processes on the Real Line}:
\textit{A Dynamic Approach}. Springer, New York.
\MR{1353912}


\bibitem{LM}
\textsc{Lions, J.-L.} and \textsc{Magenes, E.} (1968). \emph{Probl\`{e}mes aux Limites Non
Homog\`{e}nes et Applications}. Dunod, Paris.

\bibitem{LSII}
\textsc{Liptser, R. S.} and \textsc{Shiryaev, A. N.} (2000). \emph{Statistics of Random
Processes II. Applications}. Springer, New York.
\MR{1800858}

\bibitem{LSMar}
\textsc{Liptser, R. S.} and \textsc{Shiryayev, A. N.} (1989). \emph{Theory of
Martingales.} Kluwer Acad. Publ., Dordrecht.
\MR{1022664}

\bibitem{MM}
\textsc{Malliavin, P.} and \textsc{Mancino, M. E.} (2002). Fourier series method for
measurement of multivariate volatilities. \emph{Finance and Stochastics}
\textbf{6} 49--62.
\MR{1885583}

\bibitem{MP}
\textsc{Mikulevicius, R.} and \textsc{Pragarauskas, H.} (1992). On the Cauchy problem for
certain integro-differential operators in Sobolev and H\"{o}lder spaces.
\textit{Lithuanian Math.~J.} \textbf{32} \mbox{238--263.}
\MR{1246036}

\bibitem{RZ}
\textsc{Rogers, L. C. G.} and \textsc{Zane, O.} (1998). Designing and estimating models of
high-frequency data. Preprint.

\bibitem{Roz1}
\textsc{Rozovskii, B. L.} (1990). \textit{Stochastic Evolution Systems. Linear Theory and
Applications to Non-Linear Filtering}. Kluwer Acad. Publ., Dordrecht.
\MR{1135324}

\bibitem{Roz2}
\textsc{Rozovskii, B.L.} (1991). A simple proof of uniqueness for Kushner
and Zakai equations. In \textit{Stochastic Analysis} (E. Mayer-Wolf et
al., eds.) 449--458. Academic Press, Boston.
\MR{1119843}\vadjust{\goodbreak}

\bibitem{R}
\textsc{Runggaldier, W. J.} (2004). Estimation via stochastic filtering in financial
market models. In \textit{Mathematics of Finance} (G.~Yin and Q.~Zhang, eds.)
309--318. Amer. Math. Soc.,
Providence, RI.
\MR{2076550}
\end{thebibliography}
\end{document}